\documentclass[twoside]{article}

\usepackage{amsmath}
\usepackage{amsfonts}
\usepackage{amssymb}
\usepackage{color}
\usepackage[pdftex]{graphicx}
\usepackage{float}
\usepackage[english]{babel}
\usepackage[utf8x]{inputenc}
\usepackage{ucs}
\usepackage{amsmath}
\usepackage{amsfonts}
\usepackage{amssymb}
\usepackage{array}
\usepackage{listings}
\usepackage{authblk}
\usepackage[sc]{mathpazo} 
\usepackage[T1]{fontenc} 
\usepackage{microtype} 
\usepackage[margin=3.5cm,hmarginratio=1:1,top=32mm,columnsep=20pt]{geometry} 
\usepackage{multicol} 
\usepackage[hang, small,labelfont=bf,up,textfont=it,up]{caption} 
\usepackage{booktabs} 
\usepackage{hyperref} 
\usepackage{lettrine} 
\usepackage{paralist} 
\usepackage{abstract} 
\usepackage{titlesec} 
\usepackage{fancyhdr} 

\newtheorem{rem}{Remark}[section]

\newtheorem{cor}{Corollary}[section]

\newtheorem{theo}{Theorem}[section]

\newtheorem{conj}{Conjecture}[section]

\renewcommand\thesection{\Roman{section}} 
\renewcommand\thesubsection{\Roman{subsection}} 

\makeatletter
\newcommand{\@giventhatstar}[2]{#1\;\middle|\;#2}
\newcommand{\@giventhatnostar}[3][]{#1#2\;#1|\;#3#1}
\newcommand{\giventhat}{\@ifstar\@giventhatstar\@giventhatnostar}
\makeatother

\titleformat{\section}[block]{\large\scshape\centering}{\thesection.}{1em}{} 
\titleformat{\subsection}[block]{\large}{\thesubsection.}{1em}{} 
\title{\vspace{-15mm}\fontsize{24pt}{10pt}\selectfont\textbf{Extreme change-point detection}} 
\date{}

\author[1,2,3]{Kevin Bleakley}
\affil[1]{LMO, Orsay}
\affil[2]{Inria}
\affil[3]{CNRS}

\begin{document}

\maketitle

\begin{abstract}
\noindent We examine rules for predicting whether a point in $\mathbb{R}$ generated from a 50--50 mixture of two different probability distributions came from one distribution or the other, given limited (or no) information on the two distributions, and---as clues---one point generated randomly from each of the two distributions. We prove that nearest-neighbor prediction does better than chance when we know the two distributions are Gaussian densities without knowing their parameter values. We conjecture that this result holds for general probability distributions and---furthermore---that the nearest-neighbor rule is optimal in this setting, i.e., no other rule can do better than it if we do not know the distributions or do not know their parameters, or both.  
\end{abstract}

\section{Introduction}

\textbf{Warning: The conjectures in this preprint are no longer open. They were resolved in the preprint: \url{https://arxiv.org/abs/2609.20897}. New general results are also included in the new preprint.}

This work originated in trying to understand---in the most simple setting possible---what detecting one change-point in a time series truly means. Suppose we know that there is one (and exactly one) change-point in a real-valued time series of length $n$, located between the $k$th and $k+1$th data points; i.e., the first $k$ points are randomly drawn from some probability distribution $f_X$, and the last $n-k$ points from another, different, probability distribution $f_Z$. The minimal interesting setting is when $n=3$. In it, we know that the first point $X$ is generated from $f_X$, the third point $Z$ from $f_Z$, and it remains to try to work out which of the two distributions the middle point $Y$ came from (which is equivalent to predicting the change-point location). 

If we are told that the middle point is more likely to have been drawn from one of the two distributions, we can already obtain a decision rule that is better than chance: \emph{Always predict the more likely distribution, no matter the values of the three points}. 
However, with real-world data, we are unlikely to have the slightest clue about which distribution the middle point comes from. In this case, it makes sense that in the absence of prior knowledge, we treat the middle point as if it were generated from a 50--50 mixture distribution of the two distributions: $Y \sim \frac{1}{2} f_X + \frac{1}{2} f_Z$.



Note that in this setting, if we had full knowledge of the distributions $f_X$ and $f_Z$, knowing the values of the random draws $X=x$ and $Z=z$ would not in fact provide us with additional information as to whether $Y=y$ was generated from $f_X$ or $f_Z$; indeed, in this case we already have all the information we need to calculate the (optimal) Bayes classifier, (given $Y=y$), that is, the classifier that minimizes the probability of incorrect prediction of the distribution $y$ came from.
Again however, in the real world with real data, we are quite unlikely to know the distributions $f_X$ and $f_Z$ a priori. That said, in what follows, we only deal with the $n=3$ case, and these results do not generalize to $n>3$ or $\mathbb{R}^d$ for $d>1$. Nevertheless, the results are interesting in their own right, given that they seem to be highlighting curious connections between probability distributions and distances between points from these distributions.

We are interested in two distinct issues: (1) proving that there exists some decision rule that is better than a coin flip in settings where we have little or no information on the two probability distributions, and (2) proving that this decision rule is optimal, i.e., that no other prediction rule has a lower probability of incorrect prediction. In the following, we have partial answers to (1), but (2) remains completely unanswered. It may seem intuitive to some that a nearest-neighbor rule cannot be beaten in $\mathbb{R}$, but it is another thing to prove it. 



\section{A general conjecture}

Let us set the scene with a general conjecture. 

%
%

\begin{conj}
Suppose that $x$ and $z$ are drawn from arbitrary unknown probability distributions $f_X$ and $f_Z$ which are different in the sense that
$$\int_{\mathbb{R}}| f_X(w) - f_Z(w)| dw  > 0 .$$
(In some other sense could be possible too.)
Suppose also that $y$ is drawn from the mixture model
$Y \sim \frac{1}{2} f_X + \frac{1}{2} f_Z$. Let the decision rule for deciding whether $y$ was drawn from $f_X$ or $f_Y$ be the nearest neighbor rule, i.e., predict that $y$ is drawn from the distribution $f_X$ if $| x - y | < | z - y |$, and vice versa. 
Then (i) this rule is correct more than half the time, i.e., better than a coin flip, and (2) this rule minimizes the classification error, i.e., no other rule does better.
\end{conj}

\begin{rem}
This conjecture, if true, means that knowing only that $f_X$ and $f_Z$ are different (in the above sense), we can correctly predict the distribution $y$ came from more than half the time.
\end{rem}


This conjecture turns out to be non-trivial, even if we suppose that the two unknown distributions are also (unknown) probability densities. Hence, in order to take tentative steps forward, we will initially relax the problem to settings where we at least know that the two distributions are probability densities, even if we don't know their precise forms (parameters, etc.).  

\section{The Gaussian setting with equal variance}

Suppose that we are lucky and know that $f_X$ and $f_Z$ are Gaussian distributions with the same variance, but have no information on the value of this variance, nor on that of the two means, except what little can be gleaned from the three numbers $x$, $y$, and $z$. Formally, suppose that $X \sim \mathcal{N}(\mu_X,\sigma^2)$ and $Z \sim \mathcal{N}(\mu_Z,\sigma^2)$. 
In this Gaussian setting, to simplify notation and proofs, we shall write $\phi_{m,\sigma^2}$ to mean the Gaussian probability distribution function (pdf) with mean $m$ and variance $\sigma^2$, and $\Phi_{m,\sigma^2}$ the corresponding Gaussian cumulative distribution function (cdf). 

Given $x$, $y$, and $z$ drawn as before, we can ask questions like, ``Is there a decision rule for predicting $\phi_{\mu_X,\sigma^2}$ or $\phi_{\mu_Z,\sigma^2}$ for $y$ that is better than flipping a coin?'', ``Is there an optimal decision rule given what we know and don't know about the distributions and given $x$, $y$, and $z$?", and, ``How close can we get to the optimal decision rule that would be known if we knew the true values of $\mu_X$, $\mu_Z$, and $\sigma^2$?'' 

\subsection{The Bayes classifier in this setting}

If we do know $\mu_X$, $\mu_Z$, and $\sigma^2$, the Bayes classifier $C^{Bayes}$, i.e., the rule that minimizes the probability of error when predicting $\phi_{\mu_X,\sigma^2}$ vs $\phi_{\mu_Z,\sigma^2}$, is easy to calculate:  
\begin{equation}
C^{Bayes}(y) = 
\begin{cases}
\phi_{\mu_X,\sigma^2} & \text{if}\,\,\, \phi_{\mu_X,\sigma^2}(y) > \phi_{\mu_Z,\sigma^2}(y)  \\ 
\phi_{\mu_Z,\sigma^2} & \text{if}\,\,\, \phi_{\mu_X,\sigma^2}(y) < \phi_{\mu_Z,\sigma^2}(y),
\end{cases}
\end{equation}
i.e., predict the distribution of $y$ as that which has the highest density at $y$.

\begin{rem}
Note that when both Gaussian distributions have the same variance, as is the case here, the Bayes classifier is equivalent to the rule:
\begin{equation}
C^{\mu}(y) := 
\begin{cases}
\phi_{\mu_X,\sigma^2} & \text{if}\,\,\, | y - \mu_X | <  | y - \mu_Z | \\ 
\phi_{\mu_Z,\sigma^2} & \text{if}\,\,\,| y - \mu_X | >  | y - \mu_Z |,
\end{cases}
\end{equation}
where $| \cdot |$ is the absolute value.
i.e., predict the distribution whose mean is closest to $y$. Note that though this is not exactly a nearest-neighbor rule, since $\mu_X$ and $\mu_Y$ are not data points, it is not too far off one. 
\end{rem}

\subsection{Better than a coin flip}
It turns out that knowing only that we have two Gaussian distributions with the same variance, and given $x$, $y$, and $z$ drawn as before, there exist rules that are better than a coin flip if the (unknown) means $\mu_X$ and $\mu_Z$ are different. (If the two means are the same, then the two Gaussian distributions are identical and a coin flip is indeed the optimal rule.) 

\begin{theo}\label{theo1}
Suppose that $X \sim \phi_{\mu_X,\sigma^2}$, $Z \sim \phi_{\mu_X+\epsilon,\sigma^2}$, and that $Y$ is a 50--50 mixture of the two, where $\mu_X \in \mathbb{R}$, 
$\epsilon \neq 0$, and $\sigma^2$, are all unknown. Suppose we have $x$, $z$, and $y$ generated respectively from these three distributions. Then the decision rule,
\begin{equation}
C^{dist}(y) := 
\begin{cases}\label{eqth1}
\phi_{\mu_X,\sigma^2} & \text{if}\,\,\, | x - y | <  | z - y | \\ 
\phi_{\mu_X+\epsilon,\sigma^2} & \text{if}\,\,\,| x - y | >  | z - y |,
\end{cases}
\end{equation}  
has a probability of being correct greater than $1/2$.
\end{theo}

\begin{rem}
We see that it is as if we are roughly estimating $\mu_X$ by $x$ and $\mu_X + \epsilon$ by $z$ and then using these point estimates of the means in the Bayes classifier. 
\end{rem}

\noindent \textbf{Proof}. This decision rule will be correct when either $y$ is closest to $x$ and was drawn from $f_X$, or closest to $z$ and was drawn from $f_Z$. Thus we want to prove that:
\begin{equation}\label{longone}
\frac{1}{2}\mathbb{P}(\giventhat{| X - Y | < | Z - Y |}{Y \sim \phi_{\mu_X,\sigma^2}}) + \frac{1}{2}\mathbb{P}(\giventhat{| X - Y | > | Z - Y |}{Y \sim \phi_{\mu_X+\epsilon,\sigma^2}}) > \frac{1}{2}. 
\end{equation}
Since the variance of $\phi_{\mu_X,\sigma^2}$ and $\phi_{\mu_X+\epsilon,\sigma^2}$ is the same here, by symmetry it suffices to prove that
$$ P^* := \mathbb{P}(\giventhat{| X - Y | < | Z - Y |}{Y \sim \phi_{\mu_X,\sigma^2}}) > 1/2.$$
Since the value of $\mu_X$ has no influence on the following calculations, we set $\mu_X = 0$ to simplify notation. Let us define the following function of $x$ and $z$:
\begin{equation}
f(x,z) = \mathbb{P}(|x-Y| < |z-Y| \, | Y \sim \phi_{0,\sigma^2}).
\end{equation}
The values of the random variable $Y$ for which $|x-Y| < |z-Y|$ are those below $(x+z)/2$ if $x < z$, or those above $(x+z)/2$ if $x > z$. The measure of the set for which $|x-Y| < |z-Y|$ is therefore:
\begin{equation}
f(x,z) = \Phi_{0,\sigma^2}\left(\frac{x+z}{2} \right)\mathbb{1}_{\{x < z\}} +(1 -  \Phi_{0,\sigma^2})\left(\frac{x+z}{2} \right)\mathbb{1}_{\{x > z\}} .
\end{equation}
Integrating $f(x,z)$ over $x$ and $z$ will give us the probability we are looking for:
\begin{flalign*}
P^* & := P_1^* + P_2^* =
\int_{x=-\infty}^\infty \int_{z=x}^\infty  \Phi_{0,\sigma^2}\left(\frac{x+z}{2} \right)
\phi_{\epsilon,\sigma^2}(z)\phi_{0,\sigma^2}(x) dz dx & \\ 
& \qquad \qquad \qquad \qquad  + \int_{x=-\infty}^\infty \int_{z=-\infty}^x \left(1 -  \Phi_{0,\sigma^2}\left(\frac{x+z}{2} \right)\right)
\phi_{\epsilon,\sigma^2}(z)\phi_{0,\sigma^2}(x) dz dx .
\end{flalign*}
Remarking that
$$
\Phi_{0,\sigma^2}\left(\frac{x+z}{2} \right) = 
\Phi_{0,\sigma^2}\left(\frac{(x - \delta) +(z + \delta)}{2} \right)
$$
for all $\delta \in \mathbb{R}$, the following change of variable greatly simplifies the integration:
\begin{itemize}

\item For $P_1^*$  let $r = (x+z)/2$ and $\alpha = (z-x)/2$, i.e., $x = r - \alpha$ and 
$z = r + \alpha$. 

\item For $P_2^*$  let $r = (x+z)/2$ and $\alpha = (x-z)/2$  i.e., $x = r + \alpha$ and 
$z = r - \alpha$. 

\end{itemize}
In both cases, the absolute value of the Jacobian is 2. Let us concentrate on calculating the first double integral $P_1^*$ above; the second---$P_2^*$---involves almost identical calculations. We have that:
\begin{equation}\label{p1}
P_1^* = 2 \int_{r=-\infty}^\infty  \int_{\alpha=0}^\infty  \Phi_{0,\sigma^2}(r)
\frac{1}{\sqrt{2\pi}\sigma}e^{-\frac{1}{2\sigma^2}(r + \alpha - \epsilon)^2  }
\frac{1}{\sqrt{2\pi}\sigma}e^{-\frac{1}{2\sigma^2}(r - \alpha)^2  } d\alpha dr .
\end{equation}
We consider the two Gaussian pdfs as if they were functions of $\alpha$ and use the fact that 
the product of two Gaussian pdfs is proportional to another Gaussian pdf to rewrite the double integral as:
\begin{equation*}
P_1^* = \int_{r=-\infty}^\infty  \Phi_{0,\sigma^2}(r)
\frac{1}{\sqrt{2\pi}\sigma}e^{-\frac{1}{2\sigma'^2}(r - \frac{\epsilon}{2})^2  } \left(
\int_{\alpha=0}^\infty
\frac{1}{\sqrt{2\pi}\sigma'}e^{-\frac{1}{2\sigma'^2}(\alpha - \frac{\epsilon}{2})^2  } d\alpha \right) dr ,
\end{equation*}
where $\sigma`^2 =  \sigma^2/2$. The integral in $\alpha$ is nothing other than 
$1 - \Phi_{\epsilon/2,\,\sigma'^2}(0)$, or equivalently, 
$\Phi_{-\epsilon/2,\,\sigma^2/2}(0)$ (which is independent of $r$).
As for the integral in $r$, recalling the well-known result:
\begin{equation*}
\int_{u=-\infty}^\infty \Phi_{0,1}\left(\frac{u - c}{\tau_1} \right)  \phi_{0,1}\left(\frac{u - b}{\tau_2} \right) du = 
\tau_2\, \Phi_{0,1}\left(\frac{b - c}{\sqrt{\tau_1^2 + \tau_2^2}} \right) ,
\end{equation*}
it  simplifies to $\Phi_{0,1}\left(\frac{\epsilon}{\sqrt{6}\,\sigma} \right)$, and thus:
$$
P_1^* = \Phi_{-\frac{\epsilon}{2},\frac{\sigma^2}{2}}(0) \Phi_{0,1}\left(\frac{\epsilon}{\sqrt{6}\,\sigma} \right) .
$$
Essentially identical calculations show that:
$$
P_2^* = \Phi_{\frac{\epsilon}{2},\frac{\sigma^2}{2}}(0) \left(1 -  \Phi_{0,1}\left(\frac{\epsilon}{\sqrt{6}\,\sigma} 
\right) \right) .
$$
To begin to conclude, note that $\Phi_{-\epsilon/2,\,\sigma^2/2}(0) =  \Phi_{0,1}(\epsilon/\sqrt{2}\,\sigma)$ and $\Phi_{\epsilon/2,\,\sigma^2/2}(0) = 1 -  \Phi_{0,1}(\epsilon/\sqrt{2}\,\sigma)$, so that:
\begin{equation}\label{PP}
P_1^* + P_2^* = \Phi_{0,1}\left(\frac{\epsilon}{\sqrt{2}\,\sigma}\right) \Phi_{0,1}\left(\frac{\epsilon}{\sqrt{6}\,\sigma}\right) + 
\left(1 - \Phi_{0,1}\left(\frac{\epsilon}{\sqrt{2}\,\sigma}\right)\right)\left(1 - \Phi_{0,1}\left(\frac{\epsilon}{\sqrt{6}\,\sigma}\right)\right) .
\end{equation}
If $\epsilon > 0$, the domain of $\phi_{0,1}$ can be cut up into three pieces: $]-\infty,\epsilon/\sqrt{6}\,\sigma]$, $]\epsilon/\sqrt{6}\,\sigma,\epsilon/\sqrt{2}\,\sigma]$, and $]\epsilon/\sqrt{2}\,\sigma,\infty[$ with respective areas under $\phi_{0,1}$ of $j$, $k$, and $\ell$. Thus $(j + k + \ell) = 1$. Hence, 
$P_1^* + P_2^* = (j+k)j + \ell(k + \ell) = j + \ell - 2j\ell$, using the fact that  $(j + k + l)^2 = 1$ and expanding and rearranging. Then $P_1^* + P_2^* > 1/2$ is equivalent to $j > 1/2$ since $\epsilon > 0$ and $\ell < 1/2$. But $j > 1/2$ also since $\epsilon > 0$. A similar argument works when $\epsilon < 0$. $\square$


\begin{rem}
It is easy to see that that as $\epsilon \rightarrow 0$, the two Gaussians become increasingly indistinguishable and this rule tends to a probability of being correct of $1/2$, from above. When $\epsilon \rightarrow \infty$, both this rule and the Bayes classifier tend to probability of being correct of $1$. One could ask whether we could at least get an approximate idea of how well we might expect to do (between $1/2$ and 1) by roughly estimating $\epsilon$ by $z-x$, but we see in Eq.~\ref{PP} that the value of $P_1^* + P_2^*$ obtained still depends very much on the (unknown) $\sigma$.
\end{rem}


\begin{rem}
One cannot help suspecting that Eq.~\ref{PP} is related to products of areas under Gaussian densities connected to where they (or suitable normalized versions of them) cross over each other. For instance, the densities $\phi_{0,\sigma^2}$ and $\phi_{\epsilon,\sigma^2}$ cross once at $x=\epsilon/2$ and their cdfs at this point are respectively equal to $\Phi_{0,1}(\epsilon/2\sigma)$ and $1 - \Phi_{0,1}(\epsilon/2\sigma)$.  
\end{rem}


Given that we have only three data points, we may ask how this nearest neighbor classifier compares to other strategies such as maximum likelihood, CUSUM, and so on. We have the following corollary.

\begin{cor}
Under the same conditions as Theorem~\ref{theo1},  decision rule in Eq.~\ref{eqth1} corresponds to the same
solution inferred from the maximum likelihood estimator, the CUSUM method, and linear and
Gaussian kernels when running kernel change-point detection \cite{arlot2019kernel}). 
\end{cor}

\noindent \textbf{Proof}. 
Without loss of generality, suppose that $x < z$ and $y < (x+z)/2$. If you run the calculations, the posterior of the maximum likelihood estimator (requiring the EM algorithm here) essentially says that it is more likely that $y$ came from the same distribution as $x$ than from the distribution of $z$, i.e., it picks the closest point to $y$ as its classification rule. As for the CUSUM method (as described in \cite{page1954continuous}), some tedious algebra  shows that the maximum absolute value of the CUSUM criterion for the two possible change-point locations is again equivalent to the rule predicting that $y$ comes from the same distribution as the point ($x$ or $z$) closest to it. As for kernel change-point detection, the linear kernel corresponds exactly to performing least-squares minimisation of a signal with piecewise constant mean and exactly one change-point, with the same variance in both constant segments. A few lines of algebra confirm that the change-point location that minimizes the sum of squared errors is again equivalent to predicting that $y$ comes from the same distribution as the point ($x$ or $z$) closest to it. Basic calculations show that the same is also true with the Gaussian kernel. $\square$

\begin{rem}
Though tempted to conjecture that the same is true in general for other kernels, one quickly finds a counter-example: For the $1$-d polynomial kernel of degree 2: \, $k_2(u,v) = (uv)^2$, if $x=1$, $y=2$, and $z=2.9$, and thus $y$ is closer to $z$ than to $x$, the minimum of the kernel change-point criterion occurs when grouping $x$ with $y$, i.e., when putting a change-point between $y$ and $z$. This means that in some sense, in the eyes of the quadratic polynomial kernel, points further apart can be ``more similar'' to each other than points closer together. 
\end{rem}

This brings us to the heart of our other question: Can we do better than the decision rule in Theorem~\ref{theo1}, given our hypotheses? Or is it optimal, and why? If we can do demonstrably better, is this new rule optimal? And if we cannot provide a better rule, how can we prove that the ``nearest neighbor'' rule is the best we can do, since we know it is not equal to the Bayes classifier? Is this finally a question of the geometry of $1$-d space? We leave this as a conjecture.

\begin{conj}
Under the hypotheses of Theorem~\ref{theo1}, the optimal rule for deciding whether $y$ came from the same distribution as $x$ or $z$ is the nearest neighbor rule defined in Theorem~\ref{theo1}.
\end{conj}

\section{The Gaussian setting with different variances}

\begin{figure}[h!]
\begin{center}
\includegraphics[height = 7.5cm]{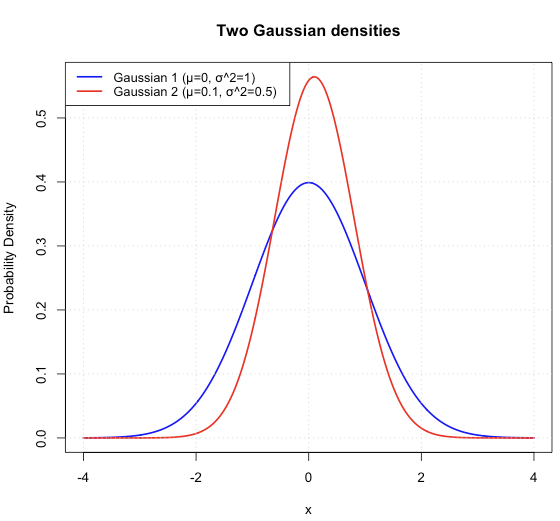}\caption{Two Gaussian densities  which do not satisfy a condition used in the proof with equal variance .}
\label{TwoGaussians}
\end{center}
\end{figure}

A natural question to ask after the previous section is whether this nearest neighbor decision rule is still valid if we know that the variances are---or could be---different for the two Gaussian distributions. Though intuition suggests that it is still valid, it turns out that certain steps in the proof of Theorem~\ref{theo1} no longer work once the two variances are different.
For instance, the symmetry argument whereby 
$$ P^* := \mathbb{P}(\giventhat{| X - Y | < | Z - Y |}{Y \sim \phi_{\mu_X,\sigma^2}}) > 1/2$$
 implies 
$$ P^* := \mathbb{P}(\giventhat{| X - Y | > | Z - Y |}{Y \sim \phi_{\mu_{X+\epsilon},\sigma^2}}) > 1/2 $$
no longer holds in general; indeed, it turns out---surprisingly---to be possible that one of these two terms can in fact be less than $1/2$ ! This occurs for example when
 $X \sim \phi_{0,1}$ and $Z \sim \phi_{0.1,0.5}$ (see Fig.~\ref{TwoGaussians}); here the nearest neighbor rule when  
$Y$ is drawn from $\phi_{0,1}$ is correct only around 44.5\% of the time! i.e., if you draw once from $\phi_{0,1}$ and once from  $\phi_{0.1,0.5}$,  a second draw from $\phi_{0,1}$ will---more than half the time---be closer to the point from the other distribution. 
However, if $\epsilon \neq 0$, then Eq.~\ref{longone} turns out to still be true in the Gaussian setting with (possibly) different variances (if one of the two probabilities is less than $1/2$, the other almost magically compensates so that their sum is greater than 1), giving us a more general result.

\begin{theo}\label{theo2}
Suppose that $X \sim \phi_{\mu_X,\sigma_X^2}$, $Z \sim \phi_{\mu_X+\epsilon,\sigma_Z^2}$, and that $Y$ is a 50--50 mixture of the two distributions, where $\mu_X \in \mathbb{R}$, 
$\epsilon \neq 0$, and $\sigma_X^2 \neq \sigma_Z^2$ are all unknown. Suppose we have $x$, $z$, and $y$ generated respectively from these three distributions. Then the decision rule,
\begin{equation}
C^{dist}(y) := 
\begin{cases}\label{eqth2}
\phi_{\mu_X,\sigma_X^2} & \text{if}\,\,\, | x - y | <  | z - y | \\ 
\phi_{\mu_X+\epsilon,\sigma_Z^2} & \text{if}\,\,\,| x - y | >  | z - y |,
\end{cases}
\end{equation}  
has a probability of being correct greater than $1/2$.
\end{theo}

\noindent \textbf{Outline of the proof}.
This proof involves two steps near the end which lack rigor; these are highlighted in the text.
Without loss of generality, rewrite $\sigma_Z^2$ as $\beta \sigma_X^2$ where $\beta$ is some unknown positive number (since  $\sigma_X^2$ and $\sigma_Z^2$ are unknown). 
As before, we first calculate
$$ P^* =  \mathbb{P}(\giventhat{| X - Y | < | Z - Y |}{Y \sim \phi_{\mu_X,\sigma_X^2}}) $$
by separating the calculation into two terms $P_1^*$ and $P_2^*$ as in Theorem~\ref{theo1}.
$P_1^*$ is now:
\begin{equation}\label{p1b}
P_1^* = 2 \int_{r=-\infty}^\infty  \int_{\alpha=0}^\infty  \Phi_{0,1}\left(\frac{r}{\sigma_X}\right)
\frac{1}{\sqrt{2\pi}\sqrt{\beta}\,\sigma_X}e^{-\frac{1}{2\beta\sigma_X^2}(r + \alpha - \epsilon)^2  }
\!\!\frac{1}{\sqrt{2\pi}\sigma_X}e^{-\frac{1}{2\sigma_X^2}(r - \alpha)^2  } d\alpha dr .
\end{equation}
We then merge the two Gaussian pdfs into a constant (w.r.t. $\alpha$) and another Gaussian pdf (a function of $\alpha$), and integrate out the latter 
as before, which gives---after copious algebra, that:
\begin{flalign*}
P_1^* & := 
2 \int_{r = -\infty}^{\infty} \Phi_{0,1}\left(\frac{r}{\sigma_X}\right)
\frac{1}{\sqrt{2\pi}\sqrt{1+\beta}\,\sigma_X}
\exp\left\{-\frac{1}{2(1+\beta)\sigma_X^2}(2r - \epsilon)^2 \right\} \times
 & \\ 
& \qquad \qquad \qquad \qquad \qquad \qquad \qquad \qquad \times \Phi_{0,1}\left(\frac{\beta-1}{\sqrt{1 + \beta}\sqrt{\beta}\,\sigma_X}r
 + \frac{1}{\sqrt{1 + \beta}\sqrt{\beta}\,\sigma_X}\epsilon \right)  dr .
\end{flalign*}
Due to a result we will shortly use, it will be useful to get the exponential term---which is more or less a Gaussian pdf---in terms of $r$ and not $2r$. After some algebra, we find that:
$$
\frac{1}{\sqrt{2\pi}\sqrt{1+\beta}\,\sigma_X}
\exp\left\{-\frac{1}{2(1+\beta)\sigma_X^2}(2r - \epsilon)^2 \right\} = 
\frac{1}{2\sigma*} \phi_{0,1}\left(\frac{r - \frac{\epsilon}{2}}{\sigma*}\right) ,
$$
where 
$$
\sigma* = \frac{1+\beta}{2}\sigma_X^2  \, .
$$
To continue, we recall a result from Owen (\cite{owen1980table}, pg.~407):
\begin{flalign}\label{yikes}
\int_{m=-\infty}^{\infty} \Phi_{0,1}(a + b m) \Phi_{0,1}(c + d  m) \phi_{0,1}(m) dm
& = \\ \nonumber
& \hspace{-5cm} \frac{1}{2}\Phi_{0,1}\left(\frac{a}{\sqrt{1 + b^2}}\right) + 
\frac{1}{2}\Phi_{0,1}\left(\frac{c}{\sqrt{1 + d^2}}\right) 
- T\left(\frac{a}{\sqrt{1 + b^2}}, \frac{c + cb^2 - abd}{a\sqrt{1 + b^2 + d^2}}\right) & \\ \nonumber
& \hspace{-3cm} - T\left(\frac{c}{\sqrt{1 + d^2}}, \frac{a + ad^2 - bcd}{c\sqrt{1 + b^2 + d^2}}\right) ,
\end{flalign}
where $T$ is Owen's $T$ function (see \cite{owen1956tables}):
$$
T(u,v) = \frac{1}{2\pi} \int_{t = 0}^{v} \frac{ e^{-\frac{1}{2} u^2}(1+t^2) }{1+t^2} dt . \qquad   (-\infty < u,\, v < +\infty) 
$$
In order to invoke this result, we make the change of variable
$m = (r - \epsilon/2)/\sigma*$ and then calculate $a$, $b$, $c$, and $d$ in our setting, 
obtaining $a = \epsilon/(2\sigma_X)$, $b = \sqrt{1+\beta}/2$, 
$c = (\sqrt{1+\beta}\,\epsilon)/(2\sqrt{\beta}\,\sigma_X)$, and 
$d = (\beta - 1)/(2\sqrt{\beta})$. Note that if $ac$ were not positive (as it clearly is in our case), there is an extra term in Eq.~\ref{yikes} which we have not shown for clarity.
Plugging these into Eq.~\ref{yikes}, we get that:
\begin{flalign}
P_1^* & = \frac{1}{2}\Phi_{0,1}\left(\frac{\epsilon}{\sqrt{5+\beta}\,\sigma_X}\right)
+ \frac{1}{2}\Phi_{0,1}\left(\frac{\epsilon}{\sqrt{1+\beta}\,\sigma_X}\right)
- T\left(\frac{\epsilon}{\sqrt{5+\beta}\,\sigma_X},\frac{3}{\sqrt{1 + 2\beta}}\right) \\
& \qquad \qquad - T\left(\frac{\epsilon}{\sqrt{1+\beta}\,\sigma_X},\frac{1}{\sqrt{1 + 2\beta}}\right) .
\end{flalign} 
One must then run through essentially the same calculations for $P_2^*$, and the result is that $P_2^*$ gives the same result as $P_1^*$ except that $\epsilon$ is replaced by $-\epsilon$. To finally calculate $P^*$ itself, we first note two things: (i) $\Phi_{0,1}(\delta) = 1 - \Phi_{0,1}(-\delta)$ for any $\delta \in \mathbb{R}$; (ii) Owen's $T$ function satisfies $T(-u,v) = T(u,v)$. Using these facts, we obtain:
$$
P^* = P_1^* + P_2^* = 1 - 2 T\left(\frac{\epsilon}{\sqrt{5+\beta}\,\sigma_X},\frac{3}{\sqrt{1 + 2\beta}}\right)
- 2 T\left(\frac{\epsilon}{\sqrt{1+\beta}\,\sigma_X},\frac{1}{\sqrt{1 + 2\beta}}\right) .
$$
Since we cannot count on symmetry here, we now have to perform this whole process again 
to calculate the other term:
$$ P^{**} =  \mathbb{P}(\giventhat{| X - Y | > | Z - Y |}{Y \sim \phi_{\mu_{X+\epsilon},\sigma_Z^2}}) .$$
Recall that $\sigma_Z^2$ is still equal to $\beta \sigma_X^2$ and that this is the \emph{same} fixed unknown $\beta$ as we have just worked with. Thus it is also true that 
$\sigma_X^2 = (1/\beta)\cdot \sigma_Z^2$. It turns out that the solution for $P^{**}$ takes exactly
the same form as that of $P^{*}$ except that $\beta$ is replaced---wherever it is found---by $1/\beta$, and $\sigma_X$ by $\sigma_Z$. We then revert the $\sigma_Z$ in $P^{**}$ back to $\sigma_X$ by multiplying by $1/\sqrt{\beta}$, and thus obtain:
$$
P^{**}  = 1 - 2 T\left(\frac{\epsilon}{\sqrt{1+5\beta}\,\sigma_X},\frac{3}{\sqrt{1 + 2/\beta}}\right)
- 2 T\left(\frac{\epsilon}{\sqrt{1+\beta}\,\sigma_X},\frac{1}{\sqrt{1 + 2/\beta}}\right) .
$$
The final probability we are looking for, which is a function of $\epsilon$ and $\beta$, is given by $P(\epsilon,\beta) = (1/2)\cdot P^{*}(\epsilon,\beta) + (1/2)\cdot P^{**}(\epsilon,\beta)$ and we must prove that this is greater than $1/2$ if $\epsilon \neq 0$ and $\sigma_X^2 \neq \sigma_Z^2$. Let us denote by $\mathcal{T}(\beta,\epsilon)$ the sum of the four Owen $T$ function integrals (ignoring the negative sign and putting aside the factor of $1/2\pi$ in front of each of them). The result will then be true if $\mathcal{T} < \pi$. Notice that each of the four Owen $T$ function integrals is of the form:
$$
\int_{0}^{f(\beta)} \frac{e^{-g(\beta) \epsilon^2 (1+t^2)}}{1+t^2} dt ,
$$
where $f$ and $g$ output only positive numbers. For any $\beta > 0$, the largest this integral can get is when $\epsilon = 0$, but since $\epsilon \neq 0$, it is more precise to say that as $\epsilon \rightarrow 0$, this integral is monotically increasing for fixed $\beta > 0$. 

If we want to bound $\mathcal{T}$ from above, we can simply bound each of the four integrals from above by setting $\epsilon = 0$ (lack of rigor \#1). Each of the four integrals then takes the more simpler form 
$$
\int_{0}^{f(\beta)} \frac{1}{1+t^2} dt ,
$$
which is in fact exactly $\text{atan}(f(\beta))$. Thus:
$$
\mathcal{T}(\beta,0) = \text{atan}\left(\frac{3}{\sqrt{1 + 2\beta}}\right)\! + \text{atan}\left(\frac{1}{\sqrt{1 + 2\beta}}\right)\! + \text{atan}\left(\frac{3}{\sqrt{1 + 2/\beta}}\right)\! + \text{atan}\left(\frac{1}{\sqrt{1 + 2/\beta}}\right). 
$$
This is a fairly nasty function to maximize analytically (lack of rigor \#2). Symbolic differentiation and root finding using Wolfram Alpha shows (or if you like, suggests) that this function $\mathcal{T}(\beta,0)$ has a unique maximum at $\beta = 1$, i.e., when 
$\sigma_X^2 = \sigma_Z^2$, which is the excluded case in the theorem's statement. Consequently, for any $\sigma_X^2 \neq \sigma_Z^2$, $\mathcal{T}(\beta,0) < \mathcal{T}(1,0) = 
2\,\text{atan}(\sqrt{3}) + 2\,\text{atan}(1/\sqrt{3}) = \pi$ by elementary properties of the atan function.
Thus $P > (1/2\pi)\cdot \pi = 1/2$ and the result is proved. $\square$

\begin{rem}
The question of the optimality of this nearest-neighbor rule under these conditions remains entirely open. 
\end{rem}


\section{More general cases}

(Below are some notes on more general cases, without proofs and potentially with  errors.)

We can now ask whether this kind of result can be extend to other densities or distributions than Gaussian ones.  For example, is this result true in general if $f_X$ and $f_Z$ have probability density functions? The following page of calculations leads to Conjecture~\ref{bigconj} below. To try and take a step in the direction of an answer, we first note that the equation to be proved (Eq.~\ref{longone}) can be rewritten:
\begin{equation}\label{longone2}
\mathbb{P}(| X - X'| < | Z - X' |) + \mathbb{P}(| X^* - Z' | > | Z^* - Z' |) \, > \, 1, 
\end{equation}
where $X, X'$, and $X^*$ are independent variables each with density $f_X$, and $Z, Z'$, and $Z^*$  independent variables each with density $f_Z$. We then remark that since
$$
\mathbb{P}(| X - X'| < | Z - X' |)  + \mathbb{P}(| X - X'| > | Z - X' |) + \mathbb{P}(| X^* - Z' | > | Z^* - Z' |) + \mathbb{P}(| X^* - Z' | < | Z^* - Z' |) = 2,$$
Eq.~\ref{longone2} will be true if and only if 
\begin{equation}\label{longone3}
\mathbb{P}(| X - X'| > | Z - X' |) + \mathbb{P}(| X^* - Z' | < | Z^* - Z' |) \, \leq \, 1.
\end{equation}
Note that the left-hand side of Eq.~\ref{longone2} can be rewritten as $\mathbb{E}[W_1]$, where
$$
W_1 = \mathbb{1}_{| X - X'| < | Z - X' |} + \mathbb{1}_{| X^* - Z' | > | Z^* - Z' |} 
$$
is a random variable that can take the values 0, 1, or 2. Similarly, the left-hand side of Eq.~\ref{longone3}
can be rewritten as $\mathbb{E}[W_2]$, where
$$
W_2 = \mathbb{1}_{| X - X'| > | Z - X' |} + \mathbb{1}_{| X^* - Z' | < | Z^* - Z' |} 
$$
is also a random variable that can take the values 0, 1, or 2. Thus the statement we wish to prove will be true if and only if $\mathbb{E}[W_2] \leq \mathbb{E}[W_1]$. 
By writing out these two expectations in the form: 
$$
\mathbb{E}[W] = 0\cdot\mathbb{P}[W = 0] + 1\cdot\mathbb{P}[W = 1] + 2\cdot\mathbb{P}[W = 2]
$$
and using independence, we quickly see that $\mathbb{E}[W_2] \leq \mathbb{E}[W_1]$ is equivalent to it being more likely that points $X'$ and $Z'$ are both closer to the other generated point from their \emph{own} distribution than both being closer to the generated point from the \emph{other} distribution. Or, to put it mathematically, $\mathbb{E}[W_2] \leq \mathbb{E}[W_1]$ if and only if:
\begin{equation}\label{ineq1}
\mathbb{P}(| X - X'| > | Z - X' |)\cdot \mathbb{P}(| X^* - Z' | < | Z^* - Z' |) \leq 
\mathbb{P}(| X - X'| < | Z - X' |)\cdot \mathbb{P}(| X^* - Z' | > | Z^* - Z' |) .
\end{equation}
Each of the four probabilities in this expression can be expanded as triple integrals and then simplified into double integrals using the same kind of steps as in the proofs of Theorems~\ref{theo1} and~\ref{theo2}. For instance,
\begin{flalign*}
\mathbb{P}(| X - X'| > | Z - X' |) & =  \int_{x=-\infty}^{\infty} \int_{z=-\infty}^{\infty} \int_{x'=-\infty}^{\infty} 
\mathbb{1}_{\{| x - x'| > | z - x' |\}} f_X(x)f_Z(z)f_X(x')df_x df_z df_{x'}  \\ 
& = \int_{x=-\infty}^{\infty} \int_{z=-\infty}^{\infty}\left[ \left(1 - F_X\left(\frac{x+z}{2}\right)\right)\mathbb{1}_{\{x < z\}} +  F_X\left(\frac{x+z}{2}\right)\mathbb{1}_{\{x > z\}}\right] f_X(x)f_Z(z) df_x df_z,  
\end{flalign*} 
where $F$ refers to the cdf of the referenced variable. If you perform this expansion for each of the four probabilities in Eq.~\ref{ineq1} and then multiply out and do some fun algebra, many terms cancel, and it turns out that the result we wish to prove overall will be true if and only if the following double integral is non-negative:
$$
\int_{r = -\infty}^{\infty} \int_{\alpha = 0}^{\infty} 
\left(F_X(r) - F_Z(r)\right)
\left(f_X(r-\alpha)f_Z(r + \alpha) - f_X(r+\alpha)f_Z(r - \alpha)  \right) 
d\alpha dr \geq 0 .
$$
Let us therefore state this as a conjecture.

\begin{conj}\label{bigconj}
Let $X \sim f_X$ and $Z \sim f_Z$ where $f_X$ and $f_Z$ are densities with cdfs $F_X$ and $F_Z$ respectively, and $X$ and $Z$ are independent.
Then:
$$
\int_{r = -\infty}^{\infty} \int_{\alpha = 0}^{\infty} 
\left(F_X(r) - F_Z(r)\right)
\left(f_X(r-\alpha)f_Z(r + \alpha) - f_X(r+\alpha)f_Z(r - \alpha)  \right) 
d\alpha dr \geq 0 .
$$
\end{conj}

\begin{rem}
Currently we have made no progress on proving or disproving this statement of the problem. One gets the feeling that if for a given $r$, $F_X(r) > F_Z(r)$ (i.e., $F_X(r) - F_Z(r)$ is a positive number), then $f_X$ has more density ``to the left'' of $r$ than $f_Z$ and consequently we could expect that---more often than not---the integral over positive $\alpha$ of $f_X(r-\alpha)f_Z(r + \alpha) - f_X(r+\alpha)f_Z(r - \alpha)$ would also be positive, and vice versa if $F_X(r) < F_Z(r)$, leading to on average a positive double integral, but this is not much of an argument. 
\end{rem}

\begin{rem}
Some possible routes to investigate this double integral: 
\begin{itemize}
\item The Wasserstein-1 distance?
\item Convolutions?
\item Rewriting the cdfs as integrals of pdfs?
\end{itemize}
\end{rem}

\begin{rem}
If this conjecture is true, the question is then whether the result of interest 
is also true for general probability distributions and not simply densities. 
\end{rem}

%
%
%

\bibliographystyle{plain}
\bibliography{biblio}

\end{document}